\documentclass[letterpaper,journal,9pt,final]{IEEEtran} 

\IEEEoverridecommandlockouts  

\usepackage[pdftex]{graphicx}
\graphicspath{{./Figs/}}
\DeclareGraphicsExtensions{.pdf}

\usepackage{cite}
\usepackage[hyperfootnotes=false]{hyperref}
\hypersetup{ 		
    pdfauthor={Florian Dorfler},        
    colorlinks=true,        
    linkcolor=black,
    urlcolor=black,
    citecolor=black     		
}

\usepackage[pdftex]{graphicx}
\graphicspath{{./fig/}}
\usepackage[usenames,dvipsnames]{color}
\usepackage{epstopdf}
\usepackage{subfigure}

\usepackage{amsmath,amssymb,mathrsfs,dsfont,mathdots}

\usepackage[nocomma]{optidef}

\newtheorem{theorem}{Theorem}[section]
\newtheorem{lemma}[theorem]{Lemma}

\newtheorem{corollary}[theorem]{Corollary}

\newtheorem{remark}[theorem]{Remark}

\newcommand\oprocendsymbol{\hbox{$\square$}}
\newcommand\oprocend{\relax\ifmmode\else\unskip\hfill\fi\oprocendsymbol}

\newcommand{\real}[0]{\mathbb R}

\DeclareSymbolFont{bbold}{U}{bbold}{m}{n}
\DeclareSymbolFontAlphabet{\mathbbold}{bbold}

\definecolor{jgreen}{rgb}{0.0, 0.5, 0.0}

\usepackage[colorinlistoftodos]{todonotes}

\usepackage[framemethod=tikz]{mdframed}

\global\mdfdefinestyle{exampledefault}{%
linecolor=gray!25,linewidth=1pt,%
backgroundcolor=gray!15, 
roundcorner=5
}

\usepackage{amsmath}

\DeclareMathOperator*{\argmin}{arg\,min}

\usepackage{accents}
\newlength{\dhatheight}

\title{On the Certainty-Equivalence Approach \\to Direct  Data-Driven LQR Design
}

\author{Florian D\"orfler, Pietro Tesi, and Claudio De Persis
\thanks{F. D\"orfler is with Department of Information Technology and Electrical Engineering, 
ETH Zurich, 8092 Zurich, Switzerland. Email: {\tt\small dorfler@ethz.ch}.
P. Tesi is with Department of Information Engineering, 
University of Florence, 50139 Florence, Italy. Email: {\tt\small pietro.tesi@unifi.it}. 
C. De Persis is with ENTEG and the J.\,C. Willems Center for
Systems and Control, University of Groningen, 8092 Groningen, 
The Netherlands. Email: {\tt\small c.de.persis@rug.nl}.
This work was supported by ETH Zurich and the SNF through the  NCCR Automation.}
}

\begin{document}

\maketitle
\thispagestyle{empty}
\pagestyle{empty}


\begin{abstract} 
The linear quadratic regulator (LQR) problem is a cornerstone of automatic control, and it has been widely studied in the data-driven setting. The various data-driven  approaches can be classified as indirect (i.e., based on an identified model) versus direct or as robust (i.e., taking uncertainty into account) versus certainty-equivalence. Here we show how to bridge these different formulations and propose a novel, direct, and regularized formulation. We start from indirect certainty-equivalence LQR, i.e., least-square identification of state-space matrices followed by a nominal model-based design, formalized as a bi-level program. We show how to transform this problem into a single-level, regularized, and direct data-driven control formulation, where the regularizer accounts for the least-square data fitting criterion. For this novel formulation we carry out a robustness and performance analysis in presence of noisy data. Our proposed direct and regularized formulation is also amenable to be further blended with a  robust-stability-promoting regularizer. In a numerical case study we compare regularizers promoting either robustness or certainty-equivalence, and we demonstrate the remarkable performance when blending both of them.
\end{abstract}



\section{Introduction}
\label{Sec: Intro}

This paper considers data-driven approaches to  {\em linear quadratic regulator} (LQR) control of linear time-invariant (LTI)  subject to process noise \cite{anderson2007optimal}. 
Data-driven control methods can be classified into {\em direct} versus {\em indirect} methods (depending 
on whether the control policy hinges upon an identified model) and {\em certainty-equivalence} 
versus {\em robust} approaches (depending on whether they take uncertainty into account) \cite{aastrom2013adaptive}. 
The relative merits of these  paradigms are well known, and we highlight the following trade-offs: For indirect methods, on the one hand, it is hard to propagate uncertainty estimates on the data through the system identification step to the control design. On the other hand, direct methods are often more sensitive to inexact data and need to be robustified at the cost of diminishing  performance. 

 For the LQR problem, a representative (though certainly not exhaustive) list of classic and recent indirect approaches (i.e., identification of a parametric model followed by model-based design) are \cite{Fiechter1997,shi2000markov,Cohen2019,mania2019certainty} in the certainty-equivalence setting and \cite{dean2019sample,umenberger2019robustLCSS,treven2020learning} in the robust case. For the direct approach we list the iterative gradient-based methods  \cite{HGGL98,fazel2018global,mohammadi2020linear}, reinforcement learning \cite{Bradtke1994}, 
behavioral methods \cite{de2019formulas}, and Riccati-based methods \cite{van2020data}
in the certainty-equivalence setting as well as \cite{de2021low,berberich2020robust,henk-ddctr-uncer} in the robust setting. 
 We remark that the world is not black and white: a multitude of approaches have successfully bridged the direct and indirect paradigms such as identification for control  \cite{hjalmarsson2005experiment,geversaa2005}, dual control \cite{feldbaum1963dual,iannelli2020structured}, control-oriented identification \cite{formentin2018core}, and regularized data-enabled predictive control \cite{bridging}. In essence, these approaches all advocate that the identification and control objectives should be blendend to regularize each other.

An emergent approach to data-driven control is borne out of the intersection of behavioral systems theory and subspace methods; see the recent survey \cite{markovskyabehavioral}. In particular, a result termed the {\em Fundamental Lemma} \cite{fundamental-lemma} implies that the behavior of an  LTI system can be characterized by the range space of a matrix containing raw time series data. This perspective gave rise to data-enabled predictive control formulations \cite{bridging,coulson2021distributionally,berberich2020data} as well as the design of explicit feedback policies \cite{de2019formulas,van2020data,de2021low,berberich2020robust}. 
Both of these are direct data-driven control approaches and robustness plays a pivotal role.

In this paper, we show how to transition between the direct and indirect as well as the robust and certainty-equivalence paradigms for the LQR problem. 
We begin our investigations with an indirect and certainty-equivalence data-driven LQR formulation posing it as model-based $\mathcal H_{2}$-optimal design, where the model is identified from noisy data by means of an ordinary least-square approach. Following \cite{bridging} we formalize this indirect approach as a bi-level optimization problem and show how to equivalently pose it as a single-level and {\em regularized} data-driven control problem. Our final problem formulation equals the one in \cite{de2019formulas} -- posing the LQR problem as a semidefinite program parameterized by data matrices -- plus an additional regularizer accounting for the least-square fitting criterion. 

The aforementioned regularizer arising from our analysis takes the form of an extra penalty term in the LQR objective function, it promotes a least-square fitting of the data akin to certainty equivalence, and it can also be interpreted as a stability-promoting term. This explains why certainty equivalence enjoys some degree of robustness to noise. With this observation and
following methods from \cite{de2021low}, we carry out a non-asymptotic analysis and give explicit conditions for robust closed-loop stability and performance bounds as a function of the \emph{signal-to-noise ratio} (SNR) for finite sample size. Different from \cite{mania2019certainty,dean2019sample}, our analysis is not restricted to Gaussian noise. In fact, we show that the certainty-equivalence approach results in stabilizing controllers whenever the SNR is sufficiently large, irrespective of the noise statistics. Further, for sufficiently large SNR, we show that the sub-optimality gap  scales linearly with the SNR. This latter result is in line with \cite{mania2019certainty,dean2019sample}, which observe that certainty equivalence performs extremely well in regimes of small uncertainty.

In a simulation case study we validate the performance of our direct, certainty-equivalence, and regularized formulation as a function of the SNR and the regularization coefficient. We also compare our formulation to a regularizer proposed in \cite{de2021low} to promote robust stability. 
The latter shows a more robust performance in case of small SNR but is inferior otherwise.
Finally, we also blend the LQR objective, our certainty-equivalence regularizer, and the robustness-promoting regularizer from \cite{de2021low} in a single direct data-driven control formulation which gives rise to a remarkable empirical performance.

The remainder of this paper is as follows. Section~\ref{sec: preliminaries} sets up the certainty-equivalence LQR problem formulation. Section~\ref{Sec: regularization} shows how to pose this problem as a direct and regularized data-driven control problem. Section~\ref{sec: analysis} presents our robustness and performance analysis. Our results are discussed in Section~\ref{sec: discussion}. Section~\ref{Sec: Numerical simulations} contains a numerical case study. Finally, Section~\ref{sec: conclusions} concludes the paper.


\section{Problem formulation: ordinary least-square identification \& certainty-equivalence LQR}
\label{sec: preliminaries}

We now formulate the model-based optimal control problem of interest and
describe the considered certainty-equivalence approach.

\subsection{Model-based linear quadratic optimal control}
\label{subsec: LQR model based}

Consider a linear time-invariant (LTI) system 
{\setlength\arraycolsep{2pt} 
\begin{eqnarray} \label{eq: ss-sys}
\left\{
\def\arraystretch{1.3}
\begin{array}{rl}
x(k+1) &= Ax(k)+Bu(k)+d(k) \\[0.2cm]
z(k) &= \left[
\begin{array}{cc}
Q^{1/2} & 0 \\ 0 & R^{1/2}
\end{array}
\right] \left[
\begin{array}{c}
x(k) \\ u(k)
\end{array}
\right] 
\end{array}
\right.\,,
\end{eqnarray}}%
where $k \in \mathbb N$,
$x \in \mathbb R^n$ is the state, $u \in \mathbb R^m$ is the control input,
 $d$ is a disturbance term, and $z$ is the performance signal of interest. We assume that $(A,B)$ is stabilizable. Finally, 
$Q \succ 0$ and $R \succ 0$ are weighting matrices. 
Here, $\succ$ ($\succeq$) and $\prec$
($\preceq$) denote positive and negative (semi)definiteness, respectively. 

The problem of interest is {\em linear quadratic regulation}  phrased as designing
a state-feedback gain $K$ that renders
$A+BK$ Schur and minimizes 
the $\mathcal H_2$-norm of the transfer function $\mathscr T(K) := d \rightarrow z$
of the closed-loop system%
\footnote{Given a stable $p \times m$ transfer function
$\mathscr T(\lambda)$ in the indeterminate $\lambda$, the $\mathcal H_2$-norm
of $\mathscr T(\lambda)$ is defined as \cite[Section~4.4]{Chen1995}:
\[
\| \mathscr T \|_2 := \sqrt{ \frac{1}{2 \pi} \int_{0}^{2 \pi} \text{trace} 
( \mathscr T (e^{j \theta})' \mathscr T (e^{j \theta}) ) d\theta }
\]}
\begin{equation} \label{eq:closed}
\left[
\begin{array}{c}
x(k+1) \\ z(k)
\end{array}
\right]  = 
\left[
\begin{array}{c|c}
A+BK & I \\ \hline 
\left[ \begin{array}{c} Q^{1/2} \\ R^{1/2} K \end{array} \right]  & 0
\end{array}
\right] \left[
\begin{array}{c}
x(k) \\ d(k)
\end{array}
\right] \,,
\end{equation}
where our notation $\mathscr T(K)$ emphasizes the dependence of the transfer function on $K$.
When $A+BK$ is Schur, it holds that \cite[Section~4.4]{Chen1995}
\begin{equation} \label{eq:trace_cost}
\| \mathscr T(K) \|_2^{2} = {\text{trace} \left( Q P + K^\top R K P \right)}\,,
\end{equation}
where $P$ is the controllability Gramian of the closed-loop system \eqref{eq:closed}, 
which coincides with the unique solution to the Lyapunov equation 
$(A+BK) P (A+BK)^\top - P + I = 0$.

The $\mathcal H_2$-norm corresponds in time domain to the energy ($\mathcal L_{2}$-norm)
of the output $z$ when impulses are applied to all input channels, and it can
be interpreted as the mean-square deviation of $z$ when 
$d$ is a white process with unit covariance, which is the classic stochastic
LQR formulation. 
Here, we view the LQR problem as a $\mathcal H_2$-optimization 
problem as our method is based on the minimization of \eqref{eq:trace_cost}. 

As shown in \cite[Section 6.4]{Chen1995}, the
controller that minimizes the $\mathcal H_2$-norm of $\mathscr T(K)$
(henceforth, \emph{optimal})
is unique and can be computed by solving a discrete-time Riccati equation \cite{anderson2007optimal}.
Alternatively, following \cite{feron1992numerical}, 
this optimal controller can be determined 
by solving the following program:%
\begin{mini}
{P\succeq I,K}{ \text{trace}\left(QP + K^{\top}RKP  \right) }{\label{eq:LQR}}{}
\addConstraint{(A+BK) P (A+BK)^{\top} - P + I}{\preceq 0 \,,}{\phantom{}}
\end{mini}
The LQR problem indeed admits many parameterizations, and the one in
\eqref{eq:LQR} can be  turned into a convex semi-definite program after a change of variables; see Section~\ref{Subsec: Tractable convex formulation} for a related  transformation.
 
We aim to compute this optimal controller in a data-driven setting when $(A,B)$
are unknown, and we have  access only to a $T$-long
stream of noisy data collected during some experiment. By \emph{noisy}
we mean that the data collected from \eqref{eq: ss-sys} are
generated with a {non-zero} disturbance $d$ that does not necessarily follow any particular statistics. 

\subsection{Subspace relations in state-space data,
ordinary least-square identification, and certainty-equivalence control}

The conventional approach to data-driven LQR is {indirect}: 
first a parametric state-space model is identified from data, and later on controllers are synthesized based on this model as in Section \ref{subsec: LQR model based}. We will briefly review this approach.
Regarding the identification task,
consider a $T$-long time series of inputs, states, 
and successor states
\begin{align*}
U_{0} &:= \begin{bmatrix} u(0) & u(1) & \dots & u(T-1)  \end{bmatrix} \in \real^{m\times T}
\,,\quad\\
D_{0} &:= \begin{bmatrix} d(0) & d(1) & \dots & d(T-1)  \end{bmatrix} \in \real^{n\times T}
\,,\quad\\
X_{0} &:= \begin{bmatrix} x(0) & x(1) & \dots & x(T-1)  \end{bmatrix}  \in \real^{n\times T}
\,,\quad\\
X_{1} &:= \begin{bmatrix} x(1) & x(2) & \dots & x(T)  \end{bmatrix}  \in \real^{n\times T}
\end{align*}
satisfying the dynamics \eqref{eq: ss-sys}, that is, 
\begin{equation}
X_{1}-D_{0} = \begin{bmatrix}B & A\end{bmatrix}  \begin{bmatrix} U_{0} \\ X_{0} \end{bmatrix} 
\label{eq: ss-sys subspace relation}\,.
\end{equation}
Let for brevity
\begin{equation*}
W_0 := \begin{bmatrix}U_{0} \\ X_{0} \end{bmatrix}\,.
\end{equation*}
We assume that the data is sufficiently rich, that is,
\begin{equation}
\text{rank} \, W_0 = n+m\,.
\label{eq: rank condition 1}
\end{equation}
The rank condition \eqref{eq: rank condition 1} is 
an identifiability condition ensuring that $(B,A)$
can be recovered from data in the noiseless case.
As shown in \cite{van2020data},  condition \eqref{eq: rank condition 1} is  generically necessary for data-driven LQR design.
In the noiseless case, this rank condition \eqref{eq: rank condition 1} is satisfied if 
the input $u$ is persistently exciting and the pair $(B,A)$ is controllable
\cite[Corollary 2]{fundamental-lemma}, thus reducing to an experiment design
condition. Condition \eqref{eq: rank condition 1} is mild also in case of noisy data,
\emph{cf.} \cite[Section 4.2]{de2021low}.

Based on $(U_{0},X_{0},X_{1})$ and under the rank condition \eqref{eq: rank condition 1},
an estimate $(\hat B, \hat A)$ of the system matrices can be obtained 
as the unique solution to the \emph{ordinary least-squares} problem%
\begin{equation}
	\begin{bmatrix} \hat B & \hat A \end{bmatrix} = 
	\argmin_{B,A} \left\| X_{1} - \begin{bmatrix}  B &  A \end{bmatrix} W_0 
	\right\|_{F} 
	= X_{1} W_0^{\dagger} ,
	\label{eq:sysid}
\end{equation}
where $\|\cdot\|_F$ denotes the Frobenius norm, and
$\dagger$ is the right inverse.

Based on the identified model in \eqref{eq:sysid}, certainty-equivalence controllers 
can be designed, \emph{i.e.}, in the LQR problem \eqref{eq:LQR}, the matrices 
$(B, A)$ are replaced by their certainty-equivalence estimates 
from \eqref{eq:sysid}. This approach can be formalized as a {\em bi-level} program:
\begin{mini}
{P \succeq I,K}{ \text{trace}\left(QP + K^{\top}RKP  \right) }{\label{eq:LQR-indirect}}{}
\addConstraint{(\hat A+\hat BK) P (\hat A+\hat BK)^{\top} - P + I \preceq 0}
\addConstraint{
\begin{bmatrix} \hat B & \hat A \end{bmatrix} = \argmin_{B,A} \left\| X_{1} - \begin{bmatrix}  B &  A \end{bmatrix} W_0 \right\|_{F}\,.
}
\end{mini}
Following the classic terminology \cite{aastrom2013adaptive}, 
we  term problem \eqref{eq:LQR-indirect} a {\em certainty-equivalence} 
and {\em indirect data-driven} control approach
and its solution $K$ a {certainty-equivalence} controller. 
It can be argued that the sequential identification-followed-by-control approach \eqref{eq:LQR-indirect} is optimal in a maximum-likelihood sense; see \cite[Section 4.2]{hjalmarsson2005experiment}. 

Note that under the identifiability condition
\eqref{eq: rank condition 1} and with noise-free data, \eqref{eq:LQR-indirect} 
is feasible and returns the optimal controller. This is because, under these
circumstances, $\hat B=B$ and $\hat A=A$ so that \eqref{eq:LQR-indirect}  
coincides with the model-based program \eqref{eq:LQR}.
In the next  sections, we present an equivalent {\em direct} data-driven control formulation and analyze its properties 
in the  case of noisy data. 

\section{Certainty equivalence as 
regularized \& direct data-driven LQR}\label{Sec: regularization}

In this section, we provide a direct data-driven formulation of certainty-equivalence LQR.
We begin our analysis by showing that the bi-level program \eqref{eq:LQR-indirect}  
can be cast as a single-level convex program with an additional regularizer accounting for implicit identification.

\subsection{Direct design \& LQR parameterization by data matrices}
\label{Subsec: data parameterization}

The approach laid out in \cite{de2019formulas} uses the subspace relations 
\eqref{eq: ss-sys subspace relation} and  \eqref{eq: rank condition 1} to 
parametrize the  LQR problem \eqref{eq:LQR} by data matrices. 
Namely, due to the rank condition 
\eqref{eq: rank condition 1}, for any  $K$, there is a matrix $G$ so that%
\begin{equation}
\begin{bmatrix} K \\ I \end{bmatrix} = W_0 G
\label{eq:cdp-pt1}\,,
\end{equation}
and due to the relation \eqref{eq: ss-sys subspace relation}
the closed-loop matrix can be parametrized directly by data matrices as%
\begin{equation}
A+BK = \begin{bmatrix} B & A \end{bmatrix}  \begin{bmatrix} K \\ I \end{bmatrix} \overset{\eqref{eq:cdp-pt1}}{=}  \begin{bmatrix} B & A \end{bmatrix} 
W_0 G \overset{\eqref{eq: ss-sys subspace relation}}{=} (X_{1}-D_{0})G
\label{eq:cdp-pt2}\,.
\end{equation}
This data-based parameterization allows us to replace 
the closed-loop matrix $A+BK$ in  \eqref{eq:LQR} by $(X_{1}-D_{0})G$ subject to the 
additional constraint \eqref{eq:cdp-pt1}. As a result, the LQR problem \eqref{eq:LQR} can be parametrized by means of the data matrices as 
\begin{mini}
{P \succeq I,K,G}{\hspace{-0.2cm} \text{trace}
\left(QP + K^{\top}RKP  \right)}{\label{eq: data-driven LQR ideal}}{}
\addConstraint{\hspace{-0.2cm}(X_{1}-D_{0})G P G^{\top}(X_{1}-D_{0})^{\top} - P + I \preceq 0 \quad}
\addConstraint{\hspace{-0.2cm}\begin{bmatrix} K \\ I \end{bmatrix} =  W_0 G}
\end{mini}
with optimal control gain $K=U_0G$. This parametrization 
is indeed a direct formulation of the LQR problem since no explicit identification
of the system matrices is involved. 
With noise-free data \eqref{eq: data-driven LQR ideal} can be efficiently implemented (after a convexification) and returns the optimal controller \cite{de2019formulas}.
With noisy data,
as $D_{0}$ is unknown, a natural approach is to disregard $D_{0}$ 
which leads to the  formulation
\begin{mini}[0]
{P \succeq I,K,G}{ \text{trace}\left(QP + K^{\top}RKP  \right)}
{\label{eq: data-driven LQR}}{}
\addConstraint{X_{1}G P G^{\top}X_{1}^{\top} - P + I \preceq 0}
\addConstraint{\begin{bmatrix} K \\ I \end{bmatrix} =  W_0 G\,,}
\end{mini}
which can be posed again as a convex program and solved efficiently; see \cite{de2021low,de2019formulas} and Section~\ref{Subsec: Tractable convex formulation} for details. Similarly to \eqref{eq:LQR-indirect},
 \eqref{eq: data-driven LQR} also enforces
some sort of certainty equivalence since the design is carried out 
as if the noise was absent.
(In \cite{de2021low}, \eqref{eq: data-driven LQR} is indeed termed 
direct certainty-equivalence approach.)  
In what follows, we show that a particular \emph{regularized} version of 
\eqref{eq: data-driven LQR} is indeed \emph{equivalent} to \eqref{eq:LQR-indirect}.

\subsection{A direct version of the certainty-equivalence LQR}
\label{Subsec: Direct Data-Driven Certainty-Equivalence LQR}
 
%
To relate \eqref{eq:LQR-indirect} and \eqref{eq: data-driven LQR}, 
consider the following program
\begin{mini}
{P \succeq I,K,G}{ \text{trace}\left(QP + K^{\top}RKP  \right) }{\label{eq:LQR-direct}}{}
\addConstraint{X_{1}G P G^{\top}X_{1}^{\top} - P + I \preceq 0}
\addConstraint{
\begin{bmatrix} K \\ I \end{bmatrix} = W_0 G
}
\addConstraint{
\left(I-W_0^{\dagger} W_0 \right) {G}  = 0\,.
}
\end{mini}
In comparison to \eqref{eq: data-driven LQR}, we have added an {\em orthogonality constraint} ensuring uniqueness of the solution $G$ in \eqref{eq:cdp-pt1}. As it will become clear from the next theorem and its corollary, problem \eqref{eq:LQR-direct} is indeed a direct version of the certainty-equivalence LQR \eqref{eq:LQR-indirect} by-passing explicit system identification 
yet robustifying the optimal control solution against noisy data -- akin to least-squares identification.\smallskip%

\begin{theorem} ({\bf\em Constraint reduction}) \label{theorem: constraint reduction} \em
Consider the direct and indirect data-driven LQR formulations \eqref{eq:LQR-indirect} and \eqref{eq:LQR-direct}, respectively. In \eqref{eq:LQR-indirect} the variables $(\hat A,\hat B)$ are uniquely determined and can be readily eliminated. Likewise, in \eqref{eq:LQR-direct} the variable $G$ is uniquely determined and can be readily eliminated. In either case, both eliminations give rise to the  identical formulation 
\begin{equation}
  \begin{aligned}
    & \underset{\textstyle P\succeq I, K}{\text{minimize}} \quad   \text{trace}\left(QP + K^{\top}RKP  \right) \\
    & \text{subject to} \\
    & \left(X_{1} W_0^{\dagger}\begin{bmatrix} K \\ I \end{bmatrix}\right) P \left(X_{1} W_0^{\dagger}\begin{bmatrix} K \\ I \end{bmatrix}\right)^{\top} - P + I \preceq 0  \,.
  \end{aligned}
  \label{eq:LQR-compact}
\end{equation}
\end{theorem}
\medskip

\emph{Proof}.
Consider the indirect data-driven problem \eqref{eq:LQR-indirect}. 
From the least-squares solution \eqref{eq:sysid} we have that
$\begin{bmatrix} \hat B & \hat A \end{bmatrix}  = X_{1}  W_0^{\dagger}$ and thus
\begin{equation*}
\hat A + \hat B K = \begin{bmatrix} \hat B & \hat A \end{bmatrix} \begin{bmatrix} K \\ I \end{bmatrix}
= X_{1} W_0^{\dagger}  \begin{bmatrix} K \\ I \end{bmatrix} 
\,.
\end{equation*}
A substitution of $\hat A + \hat B K$ in \eqref{eq:LQR-indirect} by the above formulation 
gives rise to the compact formulation \eqref{eq:LQR-compact}.

Likewise, for problem \eqref{eq:LQR-direct}, due to the orthogonality constraint 
\[
\left(I-W_0^{\dagger} 
W_0 \right) {G} = 0\,,
\]
we have that $G  \in \text{image} \,W_0^{\dagger}$. Additionally, $G$ satisfies 
$\left[\begin{smallmatrix} K \\ I \end{smallmatrix}\right] = W_0 G$, and 
$W_0$ admits a right inverse. Hence, we have $G = W_0^{\dagger}\left[\begin{smallmatrix} K \\ I \end{smallmatrix}\right]$. We are left with the compact formulation \eqref{eq:LQR-compact}.
\quad $\blacksquare$
\smallskip

\begin{corollary} ({\bf\em Equivalence of direct and indirect data-driven LQR formulations})\label{corollary: lqr equivalence} \em
Consider the direct and indirect data-driven LQR formulations \eqref{eq:LQR-indirect} and \eqref{eq:LQR-direct}, respectively. The two formulations are equivalent in the sense that the cost functions coincide and the feasible sets coincide. \quad $\blacksquare$
\end{corollary}
\smallskip

Theorem \ref{theorem: constraint reduction} suggests a transformation between the feasible sets. It is possible that the feasible sets are empty; e.g., if $(\hat B,\hat A)$ is not stabilizable. Theorem \ref{theorem: constraint reduction} remains valid though.

We term the orthogonal projector on the nullspace of $W_0$ as
\begin{equation*}
\Pi := I-W_0^{\dagger} W_0 \,.
\end{equation*}
Then,
by lifting the orthogonality constraint $\Pi G=0$ in \eqref{eq:LQR-direct} 
to the objective function, we finally 
arrive at a regularized direct data-driven LQR formulation 
mirroring that in \cite[Theorem 4.6]{bridging}:%
\begin{mini}
{P \succeq I,K,G}{ \text{trace}\left(QP + K^{\top}RKP  \right) + \lambda  \cdot \left\| \Pi {G} \right\| }{\label{eq:LQR-direct-regularized}}{}
\addConstraint{X_{1}G P G^{\top}X_{1}^{\top} - P + I \preceq 0}
\addConstraint{
\begin{bmatrix} K \\ I \end{bmatrix} = W_0 G
}
\end{mini}
where $\|\cdot\|$ is any matrix norm. 

\smallskip
\begin{theorem} ({\bf\em Regularized direct data-driven LQR})
\label{theorem: lqr regularized}\em
Consider the direct data-driven LQR formulation \eqref{eq:LQR-direct} and its regularized version \eqref{eq:LQR-direct-regularized} with parameter $\lambda \geq 0$. The two problems coincide for $\lambda$ sufficiently large. 
Otherwise, for general $\lambda \geq 0$, problem \eqref{eq:LQR-direct-regularized} lower-bounds \eqref{eq:LQR-direct}. 
\end{theorem}
\smallskip

\emph{Proof.}
The constraint $ \left\| \Pi {G} \right\|=0$ measures the distance of $G$ to the range space of $W_0$. For such a distance constraint, the equivalence of \eqref{eq:LQR-direct} and \eqref{eq:LQR-direct-regularized} for $\lambda > \lambda^{\star}$ sufficiently large is due to an exact penalization result by Clarke \cite[Proposition 2.4.3]{clarke1990optimization}. In this case, a lower bound for $\lambda^{\star}$ is the Lipschitz constant of the  objective. The latter is finite,  e.g., when reformulating \eqref{eq:LQR-direct}
as a convex program in epigraph form; see \eqref{eq:LQR-direct-convex-2} in Section~\ref{Subsec: Tractable convex formulation}.
\footnote{Alternatively, by reformulating \eqref{eq:LQR-direct} and 
\eqref{eq:LQR-direct-regularized}  as convex problems 
and certifying Slater's condition, we can leverage strong duality to show the equivalence. In this case, a lower bound for $\lambda$ is the Lagrange multiplier of the orthogonality constraint $\left\| \Pi {G} \right\|=0$. The multiplier is finite if and only if the Mangasarian-Fromovitz constraint qualification holds  \cite{gauvin1977necessary}.} For a general $\lambda \geq 0$, (not necessarily larger 
than $\lambda^{\star}$), \eqref{eq:LQR-direct-regularized} then lower-bounds 
\eqref{eq:LQR-direct}.
\quad $\blacksquare$
\smallskip

It can also be shown that \eqref{eq:LQR-direct} and  
\eqref{eq:LQR-direct-regularized} coincide for \emph{every} $\lambda \geq 0$ in the 
 case of noise-free data. We do not further elaborate on this point and 
proceed to discuss the implications of Theorem \ref{theorem: constraint reduction}
and \ref{theorem: lqr regularized}.
\smallskip

\begin{remark}[Comparison of formulations]
The standard indirect certainty-equivalence LQR problem is formulated as
the \emph{bi-level} problem \eqref{eq:LQR-indirect} consisting of sequential identification and model-based LQR. Theorem \ref{theorem: lqr regularized} shows that \eqref{eq:LQR-indirect} is equivalent to the \emph{single-level and multi-criteria} problem \eqref{eq:LQR-direct-regularized} simultaneously accounting for identification and control objectives.
This formulation is interesting in its own right, and we further elaborate on
it in Section~\ref{sec: discussion}. 
Given the equivalence of the formulations \eqref{eq:LQR-indirect}, 
\eqref{eq:LQR-direct}, \eqref{eq:LQR-compact} or \eqref{eq:LQR-direct-regularized}, the latter for $\lambda$ sufficiently large, one may wonder which is the preferred one. For now we remark that they all display similar computational performance when posed as convex programs (see Section~\ref{Subsec: Tractable convex formulation}) and defer a more in-depth discussion to Section~\ref{sec: discussion} after analyzing robustness and performance properties of certainty-equivalence LQR.
\oprocend
\end{remark}

\subsection{Tractable convex problem formulation}
\label{Subsec: Tractable convex formulation}

We briefly discuss how to convexify problem \eqref{eq:LQR-direct}
and its regularized version \eqref{eq:LQR-direct-regularized}
based on results laid out in \cite{de2019formulas}. 

First, we consider \eqref{eq:LQR-direct}. By eliminating $K = U_{0}G$ and by 
substituting $Y = GP$, we obtain that \eqref{eq:LQR-direct} is equivalent to%
\begin{mini}
{P \succeq I,X,Y}{ \text{trace}\left(QP + X \right) }{\label{eq:LQR-direct-2}}{}
\addConstraint{X_{1}G P G^{\top}X_{1}^{\top} - P + I \preceq 0}
\addConstraint{X-R^{1/2} U_{0}Y P^{-1} Y^{\top} U_{0}^{\top} R^{1/2} \succeq 0}
\addConstraint{
P  = X_0 Y
}
\addConstraint{
\Pi {Y} = 0\,,
}
\end{mini} 
with optimal controller $K = U_{0}Y P^{-1}$. By exploiting the relation $X_{0}Y=P$, 
and by applying a Schur complement, we finally arrive at the {convex} formulation
of  \eqref{eq:LQR-direct}:
\begin{mini}
{X,Y}{ \text{trace}\left(QX_{0}Y + X  \right) 
 }{\label{eq:LQR-direct-convex}}{}
\addConstraint{
\begin{bmatrix}
X_{0}Y-I & X_{1}Y \\ \star & X_{0}Y
\end{bmatrix}\succeq 0 \,,
}
\addConstraint{
\begin{bmatrix}
X & R^{1/2} U_{0} Y \\ \star & X_{0} Y 
\end{bmatrix}\succeq 0
}
\addConstraint{
\Pi {Y} = 0\,,
}
\end{mini}
with optimal controller $K = U_{0}Y (X_{0}Y)^{-1}$. 
Further, an epigraph formulation leads to the following formulation: 
\begin{mini}
{X,Y,t}{ t 
 }{\label{eq:LQR-direct-convex-2}}{}
\addConstraint{
\begin{bmatrix}
X_{0}Y-I & X_{1}Y \\ \star & X_{0}Y
\end{bmatrix}\succeq 0 \,,
}
\addConstraint{
\begin{bmatrix}
X & R^{1/2} U_{0} Y \\ \star & X_{0} Y 
\end{bmatrix}\succeq 0
}
\addConstraint{
\Pi {Y} = 0
}
\addConstraint{
t \geq \text{trace}\left(QX_{0}Y + X  \right)
\,.
}
\end{mini}
This formulation is now amenable to applying Clarke's  exact penalization result \cite[Proposition 2.4.3]{clarke1990optimization}.
After replacing 
the constraint $\Pi {Y}=0$ in \eqref{eq:LQR-direct-convex-2} by $\| \Pi Y\|=0$ and lifting it to the objective, we recover a convex formulation of the regularized problem \eqref{eq:LQR-direct-regularized}.
 
\section{Robustness and performance analysis of the certainty-equivalence LQR}
\label{sec: analysis}

\subsection{Preliminary considerations}

As shown in Section \ref{Sec: regularization}, the 
certainty-equivalence LQR problem
\eqref{eq:LQR-indirect} can be cast as  direct 
(non-sequential) control design via the single-level program \eqref{eq:LQR-direct-regularized} equipped with a regularizer.
We see that \eqref{eq:LQR-direct-regularized} searches for a solution
that satisfies the Lyapunov inequality
\begin{equation} \label{eq:Lyap}
{X_{1}G P G^{\top}X_{1}^{\top} - P + I \preceq 0}
\end{equation}
which amounts to regarding $X_{1}G$ as the closed-loop system matrix. In view of the exact relation
 $A+BK=(X_{1}-D_{0})G$ from \eqref{eq:cdp-pt2}, the stability constraint that 
should be met is actually
\begin{equation} \label{eq:Lyap_true}
{(X_{1}-D_{0})G P G^{\top}(X_{1}-D_{0})^{\top} - P + I \preceq 0}\,.
\end{equation}
In order for \eqref{eq:Lyap} to imply \eqref{eq:Lyap_true} it is sufficient 
that $G$ has small norm. %
This observation reveals one role of the regularizer $\lambda \cdot \|\Pi G\|$ 
that appears in the objective of \eqref{eq:LQR-direct-regularized}: it actually
penalizes solutions $G$ with large norm. 
In fact, Theorem \ref{theorem: lqr regularized} shows that for $\lambda$
sufficiently large the solution to \eqref{eq:LQR-direct-regularized} returns
$G = W_0^{\dagger}\left[\begin{smallmatrix} K \\ I \end{smallmatrix}\right]$
which is the least Frobenius norm $\|\cdot\|_{F}$ solution to~\eqref{eq:cdp-pt1}:%
\begin{mini*}[4]{G}{ \|G\|_{F} }{\label{eq:onG}}{}
\addConstraint{
\begin{bmatrix} K \\ I \end{bmatrix} = W_0 G
}\,.
\end{mini*}
These observations  strongly suggest that the 
certainty-equivalence LQR formulation \eqref{eq:LQR-indirect} -- which coincides 
with \eqref{eq:LQR-direct-regularized} for $\lambda$ sufficiently large --
must possess a certain degree of robustness to noise. 
Hereafter, we provide a rigorous analysis 
to this hypothesis.

\begin{remark}[Regularizations promoting stability]
\label{rem: Regularizations promoting stability}
The question when feasibility of \eqref{eq:Lyap} implies feasibility of \eqref{eq:Lyap_true} has also been studied extensively in \cite{de2021low} which proposed to regularize the  data-driven LQR problem \eqref{eq: data-driven LQR} with trace$(G P G^{\top})$. This regularizer accounts for the whole term $G P G^{\top}$ multiplying \eqref{eq:Lyap_true}, instead
of $G$ alone. %
We defer a detailed comparison of the regularizers to Sections~\ref{sec: discussion}-\ref{Sec: Numerical simulations}.
\oprocend
\end{remark}

\subsection{Main robustness and performance result}

The results of this section all refer to \eqref{eq:LQR-indirect}.
By Theorems \ref{theorem: constraint reduction} and \ref{theorem: lqr regularized}, 
\eqref{eq:LQR-indirect} is again equivalent to considering 
\eqref{eq:LQR-direct}, \eqref{eq:LQR-compact} or \eqref{eq:LQR-direct-regularized},
the latter for $\lambda$ sufficiently large. 
Denote by $K_\star$ the ground-truth optimal LQR control gain from
\eqref{eq:LQR}, and let $(\overline P,\overline K)$ be any optimal solution to 
\eqref{eq:LQR-indirect}. 
Define the \emph{signal-to-noise ratio} (SNR) by
\begin{equation} \label{eq:SNR}
\text{SNR} := \frac{\sigma_{min}(W_0)}{\sigma_{max}(D_0)} \,,
\end{equation}
where $\sigma_{min}$ and $\sigma_{max}$ denote the minimum and 
maximum singular value, respectively. Hence, \eqref{eq:SNR} gives a ratio
between the useful information $W_0$ and the useless information $D_0$.
If the identifiability
condition \eqref{eq: rank condition 1} is not satisfied, then
the SNR is zero consistent with the fact that the dynamics
cannot be identified. Under \eqref{eq: rank condition 1},
the SNR is instead well defined, strictly positive, and 
reads equivalently as
\begin{equation*} \label{eq:SNR2}
\text{SNR} = \frac{1}{\|D_0\|_{2} \|W_0^\dagger\|_{2}} \,,
\end{equation*}
where $\|\cdot\|_{2}$ is the induced 2-norm. 

We are now ready to state the robustness and performance properties 
of the certainty-equivalence LQR approach.
\smallskip

\begin{theorem} \label{thm:SNR} 
(\textbf{\emph{Closed-loop stability and performance
of certainty-equivalence LQR}})
Let $(U_0,X_0,X_1)$ be the dataset generated from
an experiment on system \eqref{eq: ss-sys},
and let the identifiability
condition \eqref{eq: rank condition 1} hold. Then,
for every 
$\varepsilon > 0$ there exists
a value $\nu>0$ such that, if the $\textrm{SNR}>\nu$, then
problem \eqref{eq:LQR-indirect} is feasible and its 
solution $(\overline P, \overline K)$ is such that $\overline K$ is stabilizing 
with sub-optimality gap
\begin{equation} \label{eq:perf}
\frac{\| \mathscr T(\overline K) \|_2^2 - \| \mathscr T(K_\star) \|_2^2}
{\| \mathscr T(K_\star) \|_2^2}
\leq \varepsilon  \,.
\end{equation}
\end{theorem}
\medskip
 
The next Section~\ref{Subsec: proof} presents the proof of Theorem \ref{thm:SNR}, 
and an in-depth discussion of the result is presented in Section \ref{sec: discussion}.

\subsection{Proof of Theorem \ref{thm:SNR}}
\label{Subsec: proof}


Recall that $K_\star$ is the ground-truth optimal LQR control from
\eqref{eq:LQR}. Let $P_\star$ be the controllability Gramian 
of the closed-loop system \eqref{eq:closed} with controller $K_\star$,
which coincides with the unique solution to the Lyapunov equation  
$(A+BK_\star) P_\star (A+BK_\star)^{\top} - P_\star + I = 0$. 
Let
\begin{equation} \label{eq:G_star}
G_\star := W_0^{\dagger}\begin{bmatrix} K_\star \\ I \end{bmatrix} 
\end{equation}
with $W_0^{\dagger}$ being the right inverse of $W_0$, which exists 
under condition \eqref{eq: rank condition 1}.
By definition of $G_\star$ and since $A+BK_\star=(X_1-D_0)G_\star$, the triplet
$(P_\star,K_\star,G_\star)$ is feasible
for \eqref{eq: data-driven LQR ideal} and, by definition of $P_\star$, satisfies 
$\|\mathscr{T}(K_\star)\|_2^2  = \text{trace} ( Q P_\star + K_\star^\top R K_\star P_\star)$.
When \eqref{eq:LQR-indirect} is feasible, we let $(\overline K, \overline P)$ be the optimal solution, and define 
\begin{equation} \label{eq:G_over}
\overline G := W_0^{\dagger}\begin{bmatrix} \,\overline K\, \\ I \end{bmatrix} \,.
\end{equation}
For compactness of notation, for a pair $(P,G)$ it is convenient to define the shorthands
\begin{eqnarray} \label{eq:MTP}
&& M := G P G^\top \nonumber \\
&& \Theta: = X_{1} M X_{1}^\top  - P \\
&& \Psi: = D_0 M D_0^\top - X_1 M D_0^\top - D_0 M X_1^\top. \nonumber
\end{eqnarray}
We finally let $(M_\star, \Theta_\star, \Psi_\star)$ and
$(\overline M, \overline \Theta, \overline \Psi)$ be defined as above
with respect to $(P_\star, G_\star)$ and $(\overline P, \overline G)$, respectively. 

Our analysis of \eqref{eq:LQR-indirect} rests on two auxiliary results
from \cite{de2021low}. We report full proofs in the Appendix as 
they differ from those in \cite{de2021low} due to the different design program.

\begin{lemma} \label{lem:tech3}
Suppose that the identifiability condition \eqref{eq: rank condition 1} holds,
and that \eqref{eq:LQR-indirect} is feasible. Let 
$(\overline P, \overline K)$ be the corresponding solution.
Let $\eta_1 \geq 1$ be a constant.
If
\begin{eqnarray} \label{eq:constr_1}
\overline \Psi \preceq \left( 1-\frac{1}{\eta_1} \right) I
\end{eqnarray}
then $\overline K$ ensures stability, and
\[
\| \mathscr T(\overline K) \|^2_2 \leq 
{\eta_1 \cdot \text{trace} ( Q \overline P 
+ \overline K{}^\top R \overline K \overline P)} \,.
\]
\end{lemma}
\medskip
 
\emph{Proof}. See the Appendix. \quad $\Box$
\smallskip

\begin{lemma} \label{lem:tech4}
Suppose that the identifiability condition \eqref{eq: rank condition 1} holds, 
and let $\eta_2\geq1$ be a constant.
If 
\begin{eqnarray} \label{eq:constr_2}
- \Psi_\star \preceq \left( 1-\frac{1}{\eta_2} \right) I
\end{eqnarray}
then \eqref{eq:LQR-indirect} is feasible and its
solution $(\overline P, \overline K)$
satisfies 
\[
\text{trace} (Q \overline P + \overline K{}^\top R \overline K \overline P  ) \leq  
\eta_2 \cdot \| \mathscr T(K_\star) \|_2^2\,.
\]
\end{lemma}
\medskip
 
\emph{Proof}. See the Appendix. \quad $\Box$
\smallskip

We now prove the feasibility statement in Theorem  \ref{thm:SNR}.
By Lemma \ref{lem:tech4}, it suffices to prove \eqref{eq:constr_2}. 
Rewrite $\Psi_\star$ as 
\[
\Psi_\star= -D_0 M_\star D_0^\top - (X_1-D_0) M_\star D_0^\top - D_0 M_\star (X_1-D_0)^\top,
\]
and notice that $(X_1-D_0) M_\star = (A+B K_\star) P_\star G_\star^\top$, where 
$M_\star = G_\star P_\star G_\star^\top$. 
Hence, \eqref{eq:constr_2} can be written as
\begin{eqnarray} \label{eq:20mod}
&& D_0 G_\star P_\star G_\star^\top D_0^\top + (A+B K_\star) P_\star G_\star^\top 
D_0^\top + \nonumber \\
&& \qquad \qquad D_0 G_\star P_\star (A+B K_\star)^\top \preceq \left( 1- \frac{1}{\eta_2} \right) I\,.
\qquad 
\end{eqnarray}
Recall that for a symmetric matrix $A$ it holds that $\|A\|_{2}\leq \alpha$ implies
$A \preceq \alpha I$. Thus, since 
$\|A+B K_\star\|_{2}$ and $\|P_\star\|_{2}$ are independent of the data, 
and by definition of $G_\star=W_0^{\dagger}\left[\begin{smallmatrix} 
K_\star \\ I \end{smallmatrix}\right]$, for every $\eta_2>1$
there is a sufficiently large SNR so that \eqref{eq:20mod} and thus
\eqref{eq:constr_2} are satisfied, hence 
such that \eqref{eq:LQR-indirect} is feasible.

To prove closed-loop stability and the sub-optimality gap \eqref{eq:perf}, we appeal to Lemma~\ref{lem:tech3} and consider \eqref{eq:constr_1}. 
Like for \eqref{eq:20mod},  it is simple to verify that 
\eqref{eq:constr_1} can be rewritten as
\begin{eqnarray} \label{eq:18mod}
&& -D_0 \overline G \overline P 
\overline G{}^\top D_0^\top - (A+B \overline K) \overline P 
\overline G{}^\top 
D_0^\top - \nonumber \\
&& \qquad \qquad D_0 \overline G
\overline P (A+B \overline K)^\top \preceq \left( 1- \frac{1}{\eta_1} \right) I 
\qquad
\end{eqnarray}

By Lemma \ref{lem:tech4},
$\text{trace} ( Q \overline P + \overline K{}^\top R \overline K \overline P  ) \leq  
\gamma$, having defined $\gamma:=
\eta_2 \cdot \| \mathscr T(K_\star) \|_2^2$. 
Since $\text{trace} ( \overline K{}^\top R \overline K \overline P  )\ge 0$, then 
$\text{trace} ( Q \overline P) 
\leq \gamma$. Moreover, $Q\succeq \sigma_{min}(Q) I$, where $\sigma_{min}(Q)$ is the smallest singular value of $Q\succ 0$.
This gives  $\overline P{}^{1/2}  Q \overline P{}^{1/2} \succeq \sigma_{min}(Q) \overline P$, 
hence $\text{trace} ( Q \overline P) = 
\text{trace} (\overline P{}^{1/2} Q \overline P{}^{1/2})\ge \sigma_{min}(Q)\, \text{trace} (\overline P) $. 
Finally, recall that
$\|A\|_2 \leq \|A\|_{\textrm{F}} = \sqrt{\text{trace} ( A^\top A )}$ for any real matrix $A$. Then
$\|\overline P\|_2 \le \text{trace} (\overline P)\le \gamma/ \sigma_{min}(Q) $. 
Analogously, we show that $\|\overline K\|_2^2 \leq \gamma/\sigma_{min}(R)$. 
In fact, $\text{trace} ( \overline K{}^\top R \overline K \overline P  ) \leq 
\gamma$. Since $\overline P\succeq I$ then
$\text{trace} ( \overline K{}^\top R \overline K \overline P )\ge
\text{trace} ( \overline K{}^\top R \overline K)$. Moreover, 
$R \succeq \sigma_{min}(R) I$, which implies $\text{trace} ( \overline K{}^\top R \overline K)\ge \sigma_{min}(R)\, 
\text{trace} ( \overline K{}^\top \overline K)$. Hence,
$\text{trace} ( \overline K{}^\top \overline K)\ge \|\overline K\|_2^2$, so that
$\|\overline K\|_2^2\le \gamma/\sigma_{min}(R)$. 
Thus, all the terms in \eqref{eq:18mod}
are upper bounded by data-independent quantities except for 
$D_0 \overline G= D_0 W_0^\dagger\left[\begin{smallmatrix} 
\overline K \\ I\end{smallmatrix}\right]$, hence  
\eqref{eq:18mod} is satisfied when the SNR is sufficiently large. 

We have just proved that if the SNR is sufficiently large then
\eqref{eq:constr_1} and \eqref{eq:constr_2} hold. 
By Lemma \ref{lem:tech4}, \eqref{eq:LQR-indirect} is feasible  
and its solution 
satisfies $\text{trace} ( \overline P + \overline K{}^\top R \overline K \overline P  ) 
\leq \eta_2 \cdot \| \mathscr T(K_\star) \|_2^2$. 
Further, by Lemma \ref{lem:tech3}, $\overline K$ ensures
closed-loop stability with
$\| \mathscr T(\overline K) \|_2^2 \leq 
\eta_1 \cdot \text{trace} ( \overline P + \overline K{}^\top R \overline K \overline P))$.
By combining the two inequalities we obtain 
$\| \mathscr T(\overline K) \|_2^2 \leq 
\eta_1 \eta_2 \| \mathscr T(K_\star) \|_2^2$ which, 
written in a ratio form and with the substitution $ \varepsilon = \eta_1 \eta_2-1$,
gives the result. \quad 
\smallskip

\section{Discussion} \label{sec: discussion}

\subsection{Scaling of performance as SNR$^{-1}$}

Theorem \ref{thm:SNR} provides a robustness and performance analysis of  
certainty-equivalence LQR as a function of the SNR \eqref{eq:SNR} and independent of any noise statistics. 
The result aligns well with recent work on data-driven LQR design \cite{mania2019certainty,dean2019sample}. There
it is observed that in regimes of small estimation errors
the certainty-equivalence approach performs extremely favorably, with a 
sub-optimality gap scaling linearly with the estimation error
(and quadratically as the 
uncertainty further decreases \cite{mania2019certainty}).
Theorem \ref{thm:SNR} indeed shows that for sufficiently large SNR the left-hand sides of \eqref{eq:constr_1} and \eqref{eq:constr_2}
decreases as SNR$^{-1}$ (\emph{cf.} \eqref{eq:20mod} and \eqref{eq:18mod}).
This implies that the right-hand side of \eqref{eq:perf} controling
the sub-optimality gap decreases as SNR$^{-1}$, too.

We emphasize that these conclusions are independent of any noise statistics and 
consistent with the fact that the estimation error
 from the ordinary least-square identification \eqref{eq:sysid} scales as SNR$^{-1}$:
\[
\left\| \begin{bmatrix} \hat B & \hat A \end{bmatrix} -
\begin{bmatrix}  B &  A \end{bmatrix} \right\|_2 = \| D_0 W_0^{\dagger} \|_2 \leq 
\frac{1}{\text{SNR}}
\]

\subsection{Certainty-equivalence versus robust control}

To cope with severe uncertainties (e.g., larger noise levels),
many approaches that address LQR\,/\,$\mathcal{H}_2$-control design
{explicitly} consider robustness against uncertainty  
\cite{dean2019sample,berberich2020robust,henk-ddctr-uncer,umenberger2019robustLCSS,treven2020learning,de2021low}. 
This is achieved by deriving ellipsoidal 
confidence regions describing the
set of all systems consistent with the observed data 
and priors on the noise (this can be done deterministically,
in a set-membership sense \cite{berberich2020robust,henk-ddctr-uncer},
or probabilistically \cite{dean2019sample,umenberger2019robustLCSS,
treven2020learning}) and by subsequently applying design tools from robust control,
such as the \emph{System Level Synthesis} (SLS) 
\cite{dean2019sample,treven2020learning}
or the \emph{$S$-procedure} \cite{umenberger2019robustLCSS,henk-ddctr-uncer}.
If feasible, these approaches explicitly enforcing robustness typically outperform  
certainty-equivalence design when it comes to ensuring closed-loop stability,
but they are much less performing in regimes of small uncertainty; 
\emph{e.g.}, \emph{cf.} \cite{dean2019sample} for  simulations comparing 
certainty equivalence with SLS, and \cite[Section 2.1]{mania2019certainty} 
for a theoretical comparison. Along similar lines, \cite{de2021low} proposes a regularizer promoting robust stability; see Remark~\ref{rem: Regularizations promoting stability}.

Section~\ref{Sec: Numerical simulations} is explicitly devoted to a numerical case study comparing robust and certainty-equivalence approaches.

\subsection{Sample complexity and Gaussian noise statistics}

The analysis of Theorem \ref{thm:SNR}
is non-asymptotic. In particular,    
a controller with nearly optimal performance can be synthesized even for a small sample size. 
The only constraint is due to the identifiability condition \eqref{eq: rank condition 1} 
which requires a minimum of $T\geq n+m$ samples.

As discussed above, in the noisy case the suboptimality gap scales as SNR$^{-1}$, independent of the noise statistics. If the noise follows a Gaussian distribution, then by averaging data matrices from multiple experiments, a high-confidence bound can be explicitly stated for the SNR: namely, it decays inversely proportional to the square root of the number of experiments; see  \cite[Section 6.2]{de2021low} for details.

\subsection{Comparison of direct and indirect problem formulations}

Our analysis shows that  the indirect (bi-level) and 
certainty equivalence LQR problem \eqref{eq:LQR-indirect} can be cast equivalently as the single-level problem \eqref{eq:LQR-direct} or  \eqref{eq:LQR-direct-regularized} for sufficiently large regularization coefficient. We want to briefly point out the merits of the latter single-level formulations over the conventional bi-level formulation.

First, a single-level formulation 
 leads to a robustness and performance analysis  which is arguably
simpler than the one that we obtain in an indirect and parametric model-based setting, e.g., 
through 
perturbation analysis of LMIs or Riccati equations (\emph{e.g.}, \emph{cf.} \cite{mania2019certainty}). 
Indeed, our proofs of Lemma \ref{lem:tech3} and \ref{lem:tech4}  
use the equivalence of \eqref{eq:LQR-indirect} and
\eqref{eq:LQR-direct}.
%
%
Moreover, the uncertainty quantifications in identification and control are usually incompatible since the former are often stochastic and the latter typically require robust formulations.
In contrast, our novel direct data-driven formulations \eqref{eq:LQR-direct}--\eqref{eq:LQR-direct-regularized} are amenable to a theoretic analysis in presence of noise. Further, they allow to directly map uncertainty on the data  to the control problem and lend themselves towards further robustifications, e.g., augmenting them with robustness-promoting regularizers; see Section~\ref{Sec: Numerical simulations}.

Second, the regularized formulation \eqref{eq:LQR-direct-regularized} has its own merits over hard-coding the least-squares objective as a constraint, as it is done in \eqref{eq:LQR-direct} or \eqref{eq:LQR-compact}. Namely, it permits to modify the LQR objective in 
a smooth manner. Intuitively, we can trade off performance and robustness objectives by changing the regularizer. A case study in Section~\ref{Sec: Numerical simulations} shows the remarkable performance when blending the certainty-equivalence regularizer $\|\Pi G\|$  with the stability-promoting regularizer from Remark~\ref{rem: Regularizations promoting stability}. We believe that this line of work deserves consideration beyond LQR to better understand multi-objective problems, where performance and robustness goals coexist.

Third and finally, the fact that the bi-level certainty-equivalence  LQR problem \eqref{eq:LQR-indirect}  
can be cast as the single-level multi-criteria problem \eqref{eq:LQR-direct-regularized} is interesting in its own right. Problem \eqref{eq:LQR-direct-regularized}
 simultaneously accounts for identification and control goals similar to identification for control  \cite{hjalmarsson2005experiment,geversaa2005}, dual control \cite{feldbaum1963dual,iannelli2020structured}, control-oriented identification \cite{formentin2018core}, and data-enabled predictive control \cite{bridging} all advocating that identification and control goals should regularize each other. 

\subsection{Data-dependent stability test}

Our analysis
also gives a method to certify closed-loop stability from data whenever
we know an upper bound on the noise magnitude. In fact, 
if $\|D_0\| \leq \delta$ for some known $\delta >0$, then 
(using the notation from Section~\ref{Subsec: proof}) stability can be certified via the condition
\begin{eqnarray} \label{eq:constr_1_test}
\delta^2 \|\overline M\| + 2 \delta \|X_1 \overline M\| \leq 1-\frac{1}{\eta_1} \,.
\end{eqnarray}
If fulfilled, this condition implies \eqref{eq:constr_1}, which 
guarantees closed-loop stability, in agreement with Lemma \ref{lem:tech3}. Further, when augmenting any of the formulations \eqref{eq:LQR-indirect},
\eqref{eq:LQR-direct}-\eqref{eq:LQR-direct-regularized}
with the constraint $K=0$, \eqref{eq:constr_1_test} provides a  
test for assessing open-loop stability from noisy data.

\section{Numerical simulations}
\label{Sec: Numerical simulations}

In this section, we exemplify our main theoretical findings through simulations.
Consider the system proposed  in \cite[Section 6]{dean2019sample} given by
\begin{eqnarray*}
A = \begin{bmatrix}
1.01 & 0.01 & 0 \\
0.01 & 1.01 & 0.01 \\ 
0 & 0.01 & 1.01
\end{bmatrix}, \quad
B = I\,.
\end{eqnarray*}
These dynamics correspond to a discrete-time marginally unstable Laplacian system.
As weight matrices, we select $Q=I$ and $R=10^{-3}I$. A small input weight $R$  relative to the
state weight $Q$ favors stabilizing solutions \cite[Section 5]{de2021low}. In particular, this choice makes it 
possible to find stabilizing controllers even from a single experiment.

\begin{figure}[t!]
\includegraphics[width=.5\textwidth]{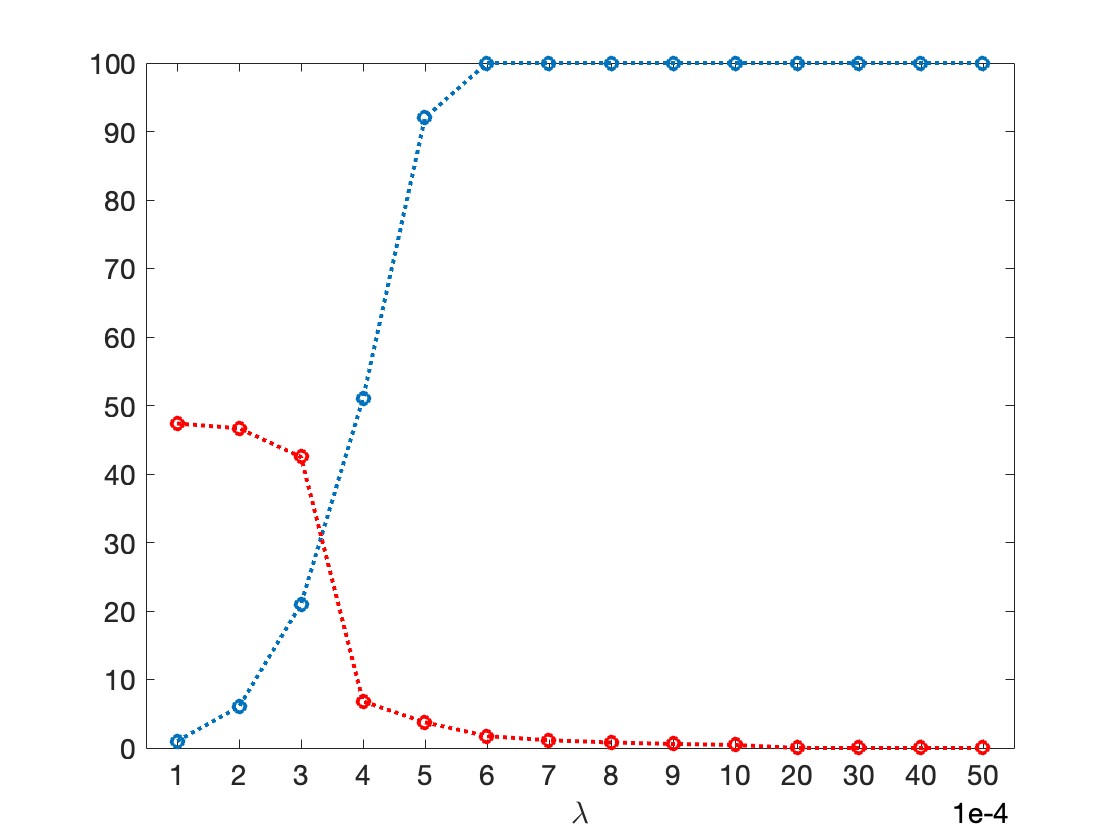}\\
\includegraphics[width=.5\textwidth]{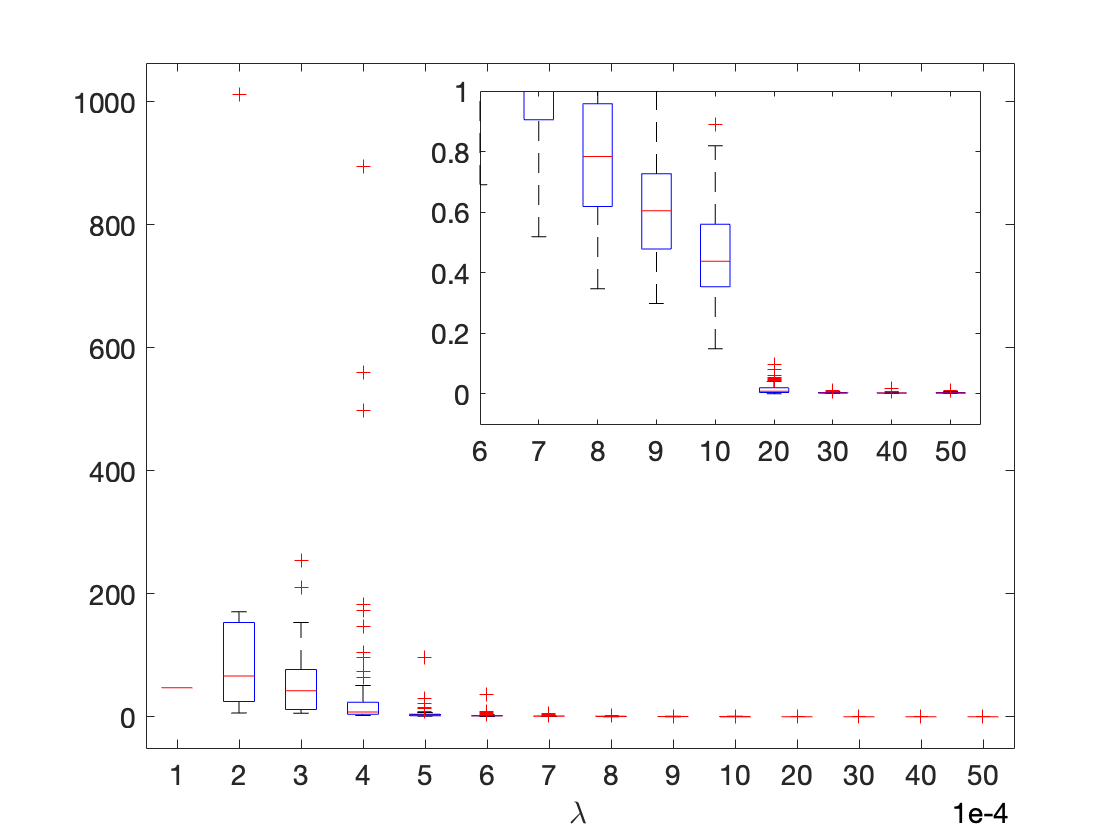} 
\caption{Performance of the direct regularized approach \eqref{eq:LQR-direct-regularized} as a function of $\lambda$.
Top panel: the blue curve the displays percentage $\mathcal S$ of stabilizing controllers,
while the red curve reports the median empirical error $\mathcal M$.
Bottom panel: box plot of the empirical errors \eqref{eq:perf_trial}.
In agreement with Theorem \ref{theorem: lqr regularized}, 
the approach \eqref{eq:LQR-direct-regularized} coincides with 
\eqref{eq:LQR-indirect} (equivalently \eqref{eq:LQR-direct})
for $\lambda$ sufficiently large. For this data set, $\lambda \geq 0.0028$ gives minimum error $\min_k \mathcal E_k = 5.922e$-$4$,
maximum error $\max_k \mathcal E_k = 0.0094$, and
median error $\mathcal M=0.0026$.
}
\label{fig:varying_lambda}
\end{figure}

Figure \ref{fig:varying_lambda} shows the results obtained with
the approach \eqref{eq:LQR-direct-regularized} as we vary the regularization coefficient $\lambda$.
We consider $100$ trials. For each trial we run an experiment on the system with 
input $u \sim \mathcal N(0,I)$ and disturbance $d \sim \mathcal N(0,0.01 I)$, 
and we collect $T=20$ state and input samples. We let $\overline K{}^{(k)}$ be the controller 
obtained in $k$-th trial. Whenever $\overline K{}^{(k)}$ is stabilizing, we define the empirical error
\begin{equation} \label{eq:perf_trial}
\mathcal E_k := \frac{\| \mathscr T(\overline K{}^{(k)}) \|_2^2 - \| \mathscr T(K_\star) \|_2^2}
{\| \mathscr T(K_\star) \|_2^2}
\end{equation}
We denote by 
$\mathcal S$ the percentage of times that we find a stabilizing controller
and by $\mathcal M$ the median of $\mathcal E_k$ through all trials.
We consider the median because it is more robust to outliers: extreme values of $\mathcal E_k$
 are due to a particular noise realization; see the box plot in  Figure~\ref{fig:varying_lambda}.

Figure \ref{fig:varying_lambda} confirms 
that regularization is indeed needed and that the regularized certainty-equivalence approach \eqref{eq:LQR-direct-regularized} 
(\emph{cf.} Theorem \ref{theorem: lqr regularized}) is robust to noisy data
and achieves excellent performance for sufficiently large $\lambda$. Namely, $\mathcal S = 100\%$ and $\mathcal M = 0.0026$ for $\lambda \geq 0.0028$. Further, the box plot evinces that  performance becomes reliable (i.e., rapidly diminishing outliers) as $\lambda$ increases. These findings are aligned with those in \cite{bridging}, and the performance of the certainty-equivalence approach is indeed remarkable considering that 
each trial involves only a single experiment with $T=20$ samples.

\begin{table*}[t] \label{sim:table2}
\begin{center}
\footnotesize
 \begin{tabular}{|c|c|c|c|c|c|c|} 
 \hline & 
 $\sigma=0.01$  & 
 $\sigma=0.1$ & 
 $\sigma=0.3$ & 
 $\sigma=0.7$ & 
 $\sigma=1$ \\ [0.5ex] 
 & 
(SNR $>15$dB) & 
(SNR $\in[5,10]$dB) & 
(SNR $\in[0,5]$dB) & 
(SNR $\approx 0$dB) & 
(SNR $<-5$dB) \\ [0.5ex] 
 \hline\hline
Certainty-equivalence &  
 $\mathcal S = 100\%$ &  
 $\mathcal S = 100\%$ &  
 $\mathcal S = 100\%$ &
 $\mathcal S = 97\%$ & 
 $\mathcal S = 84\%$  \\ [0.5ex] 
 approach \eqref{eq:LQR-indirect} & 
 $\mathcal M = 2.5599e$-$05$ &  
 $\mathcal M = 0.0026$ &  
 $\mathcal M = 0.0237$ &
 $\mathcal M = 0.1366$ & 
 $\mathcal M = 0.2596$ \\ [0.5ex] 
 \hline
 Robust approach \cite{de2021low} &  
 $\mathcal S = 100\%$ & 
 $\mathcal S = 100\%$ &  
 $\mathcal S = 100\%$ & 
 $\mathcal S = 100\%$ & 
 $\mathcal S = 100\%$ \\ [0.5ex] 
 &  
 $\mathcal M = 0.0035$ & 
 $\mathcal M = 0.0074$ &  
 $\mathcal M = 0.0369$ & 
 $\mathcal M = 0.2350$ & 
 $\mathcal M = 0.6270$ \\ [0.5ex] 
 \hline
 Mixed regularization &  
 $\mathcal S = 100\%$ & 
 $\mathcal S = 100\%$ &  
 $\mathcal S = 100\%$ & 
 $\mathcal S = 100\%$ & 
 $\mathcal S = 100\%$ \\ [0.5ex] 
combining \eqref{eq:LQR-indirect} with \cite{de2021low}&  
$\mathcal M = 0.0035$ & 
$\mathcal M = 0.0060$ &  
$\mathcal M = 0.0243$ & 
$\mathcal M = 0.1242$ & 
$\mathcal M = 0.2912$ \\ [0.5ex] 
 \hline
\end{tabular}
\smallskip\smallskip
\caption{Comparison among different approaches as we vary the noise variance $\sigma^2$.
}
\label{table2}
\end{center}
\end{table*}
\normalsize

Table \ref{table2} shows the performance of the certainty-equivalence approach 
for different values of the noise variance, i.e., for different SNR values.
Further, we compare the certainty-equivalence approach 
with the robust approach proposed in \cite{de2021low}, namely program \eqref{eq: data-driven LQR} with 
regularizer trace$(G P G^{\top})$; see Remark~\ref{rem: Regularizations promoting stability}. We refer the interested 
reader to \cite{dean2019sample} for numerical simulations comparing the
certainty-equivalence approach with the robust approach based on SLS. 
In line with the discussion of Section \ref{sec: discussion} and with the conclusions of  
\cite{dean2019sample,mania2019certainty}, the numerical simulations indicate that 
certainty-equivalence controllers are less robust but,
when stabilizing, significantly outperform robust controllers. Finally, Table \ref{table2}
also shows the remarkable performance obtained 
when blending the certainty-equivalence approach with \cite{de2021low},
namely program \eqref{eq: data-driven LQR} with the regularizer $\|\Pi G\|$ + trace$(G P G^{\top})$.
Understanding how to properly select and combine different regularizers
deserves consideration beyond LQR design.

\section{Conclusions}
\label{sec: conclusions}

We have proposed a novel, direct, and regularized data-driven LQR formulation that is equivalent to the classic indirect certainty-equivalence LQR, where a least-squares identification of the state matrices is followed by a nominal model-based design. We have formally shown this equivalence and provided a robustness and performance analysis in presence of noisy data. Our formulation is also amenable to be augmented with a robustness-promoting regularization. By varying the regularization coefficients, we can interpolate between robust and certainty-equivalence design and also recover the indirect approach for sufficiently large  coefficient. A numerical case study has illustrated the merits of the different formulations and highlighted the remarkable performance obtained with a mixed  regularization.

Surprisingly, given the remarkable empirical performance and theoretical tractability, we note that our approach is arguably {\em simple} -- both in derivation and implementation. We envision that this simplicity makes our work amenable to extensions to different system classes, identification criteria, regularization terms,  or control objectives.

\appendix

\subsection{Proof of Lemma \ref{lem:tech3}}

Suppose that \eqref{eq:LQR-indirect} is feasible and let $(\overline P,\overline K)$
be the optimal solution. Let $\overline G := W_0^{\dagger}\left[\begin{smallmatrix} \overline K \\ I \end{smallmatrix}\right]$,
as in \eqref{eq:G_over}. By \eqref{eq: rank condition 1}, we have  
$\begin{bmatrix} \hat B & \hat A \end{bmatrix}  = X_{1}  W_0^{\dagger}$.
Hence, by definition of $\overline G$,
the triplet $(\overline P,\overline K,\overline G)$ is feasible for \eqref{eq:LQR-direct}.
In particular, $(\overline P,\overline G)$ satisfies  
\begin{equation} \label{eq:app1}
X_{1} \overline G \overline P {\overline G}^\top X_{1}^\top - \overline P + I \preceq 0 \,.
\end{equation}
We will now exploit \eqref{eq:app1} to show that, under \eqref{eq:constr_1}, $(\eta_1 \overline P,\overline K,\overline G)$ is
feasible for \eqref{eq: data-driven LQR ideal}. 
To this end, rewrite \eqref{eq:app1} compactly as
$\overline \Theta + I \preceq 0$ where $\overline \Theta$ is as in 
\eqref{eq:MTP}. 
We have
\begin{eqnarray*}
&& \eta_1 \overline \Theta + \eta_1 \overline \Psi + I =  \nonumber \\
&& \eta_1 ( \overline \Theta + \overline \Psi ) + \eta_1 I + (1 - \eta_1) I  = \\
&& \eta_1 (\overline \Theta + I) + \eta_1 \overline \Psi + (1 - \eta_1) I \preceq 0 \,,
\nonumber
\end{eqnarray*}
where the inequality follows from $\eta_1 (\overline \Theta + I) \preceq 0$ and \eqref{eq:constr_1}.
Thus $\eta_1 \overline \Theta + \eta_1 \overline \Psi + I \preceq 0$,
therefore $(\eta_1 \overline P,\overline G)$ satisfies the first constraint
of \eqref{eq: data-driven LQR ideal}, namely
\begin{equation} \label{eq:app2}
(X_{1}-D_0) \overline G (\eta_1 \overline P) {\overline G}^\top (X_{1}-D_0)^\top - (\eta_1 \overline P) + I \preceq 0 \,.
\end{equation}
By definition of 
$\overline G$, the pair $(\overline K,\overline G)$
satisfies also the 
second constraint of \eqref{eq: data-driven LQR ideal}. 
The result then follows from \cite[Lemma 2]{de2021low}. \quad $\blacksquare$

\subsection{Proof of Lemma \ref{lem:tech4}}

As shown after \eqref{eq:G_star}, the condition \eqref{eq: rank condition 1} ensures that 
$(P_\star, K_\star, G_\star )$ is feasible for \eqref{eq: data-driven LQR ideal}. Thus we have 
$\Theta_\star  + \Psi_\star  + I \preceq 0$, so that
\begin{eqnarray*}
&& \eta_2 \Theta_* + I = \nonumber \\
&& \eta_2 (  \Theta_\star  + \Psi_\star  ) - 
\eta_2 \Psi_\star   + \eta_2 I + (1 - \eta_2) I = \\
&& \eta_2 (  \Theta_\star  + \Psi_\star  + I) 
- \eta_2 \Psi_\star  + (1 - \eta_2) I \preceq 0 \,,
\end{eqnarray*}
where the inequality follows from $\eta_2 (  \Theta_\star  + \Psi_\star  + I) \preceq 0$ 
and \eqref{eq:constr_2}. Thus $(\eta_2 P_\star,K_\star, G_\star)$ is 
feasible for \eqref{eq:LQR-direct}, in particular 
\begin{equation} \label{eq:app3}
X_{1} G_\star (\eta_2 P_\star) G_\star^\top X_{1}^\top - (\eta_2 P_\star) + I \preceq 0 \,.
\end{equation}
Since $G_\star = W_0^{\dagger}\left[\begin{smallmatrix}  K_\star \\ I \end{smallmatrix}\right]$
and because $\begin{bmatrix} \hat B & \hat A \end{bmatrix}  = X_{1}  W_0^{\dagger}$,
we have that $(\eta_2 P_\star,K_\star)$ is feasible for \eqref{eq:LQR-indirect}. 
The claim then follows since
$(\overline P,\overline K)$ is optimal for \eqref{eq:LQR-indirect}
and because the cost of the solution $(\eta_2 P_\star,K_\star)$ is
$\text{trace} ( Q \eta_2 P_\star  + K_\star^\top R K_\star \eta_2 P_\star) 
= \eta_2 \cdot \| \mathscr T(K_\star ) \|_2^2$. \quad $\blacksquare$

\bibliographystyle{IEEEtran}
\bibliography{refs}


\end{document}